\documentclass[11pt]{article}
\usepackage{amsmath}
\usepackage{amsthm}
\usepackage{amssymb}
\usepackage{concrete}
\def\pf{{\bf Proof. }}
\def\beq{\begin{equation}}
\def\eeq{\end{equation}}
\def\g{\gamma}

\def\al{\alpha}

\def\nd{\noindent}
\def\d{\delta}

\def\p{\partial}
\def\hf{\hfill{$\Box$}}
\def\bn{{\bf n}}
\def\im{{\rm Int}(M)}
\def\<{\leq}
\def\>{\geq}

\newtheorem{lem}{Lemma}[section]
\newtheorem{prop}{Proposition}[section]

\begin{document}
\title{\bf The Riemannian Manifolds with Boundary and Large Symmetry}
\author{\sc Zhi Chen, Yiqian Shi and Bin Xu$^\dagger$}
\date{}
\maketitle
\begin{abstract}
\noindent Sixty years ago, S. B. Myers and N. E. Steenrod ({\it
Ann. of Math.} {\bf 40} (1939), 400-416) showed that the isometry
group of a Riemannian manifold without boundary has a structure of
Lie group. Recently A. V. Bagaev and N. I. Zhukova ({\it Siberian
Math. J.} {\bf 48} (2007), 579-592) proved the same result for a
Riemannian orbifold. In this paper, we firstly show that the
isometry group of a Riemannian manifold $M$ with boundary has
dimension at most $\frac{1}{2}\,\dim\,M\,(\dim\,M-1)$. Then we
completely classify such Riemannian manifolds with boundary that
their isometry groups attain the preceding maximal dimension.
\end{abstract}


\noindent
 {\bf Key  Words and Phrases:} Riemannian manifold with boundary, isometry, rotationally symmetric metric,
principal orbit.

\footnote{\hspace{-0.5cm}$^\dagger$ Bin Xu is the correspondant
author, whose e-mail address is \tt{bxu@ustc.edu.cn} .}

%%%%%%%%%%%%%%%%%%%%%%%% INTRODUCTION %%%%%%%%%%%%%%%%%%%%%%%%%%%%%
\section{Introduction}

Let $M$ be a connected smooth Riemannian manifold with or without
boundary. A priori, there exists two definitions of isometry on
$M$. The first one is given above as we think $M$ of as a metric
space. The second is defined to be a diffeomorphism of $M$ onto
itself which preserves the metric tensor. In case of Riemannian
manifolds without boundary, these two definitions are equivalent
according to Myers and Steenrod \cite{MS} in 1939 (see also pp.
169-172 in Kobayashi-Nomizu \cite{KN} for a proof). Moreover,
Myers and Steenrod \cite{MS} proved the following result on the
isometry group of a Riemannian manifold
without boundary. \\

\nd{\bf Fact 1.1.} {\it Let $M$ be a connected smooth Riemannian
manifold without boundary.  Then the isometry group $I(M)$ is a
Lie transformation group with respect to the compact-open
topology. For each $x\in M$, the isotropy
subgroup $I_x(M)$ is compact. If $M$ is compact, $I(M)$ is also compact.}\\

S. Kobayashi (see p. 15 and p. 41 in \cite{Ko}) gave a different
proof to Fact 1.1 by the concept of $G$-structure from the
original one by Myers-Steenrod. Furthermore Kobayashi proved that
there exists a natural embedding of the isometry group $I(M)$ into
the orthonormal frame bundle $O(M)$ of $M$ such that $I(M)$
becomes a closed submanifold of $O(M)$. It is this submanifold
structure that makes $I(M)$ into a Lie transformation group.
Following this idea, he  further proved \\

\nd{\bf Fact 1.2.} (see pp. 46-47 in \cite{Ko}) {\it Let $M$ be an
$n$-dimensional connected Riemannian manifold without boundary.
Then the isometry group $I(M)$ has dimension at most $(n+1)n/2$.
If $\dim\, I(M)=n(n+1)/2$, then $M$ is isometric to one of
the following spaces of constant sectional curvature{\rm :}\\
{\rm (a)} An $n$-dimensional Euclidean space ${\bf R}^n$.\\
{\rm (b)} An $n$-dimensional unit sphere $S^n=\{x\in{\bf R}^{n+1}:|x|=1\}$ in ${\bf R}^{n+1}$.\\
{\rm (c)} An $n$-dimensional projective space ${\bf R}P^n=S^n/\{\pm 1\}$.\\
{\rm (d)} An $n$-dimensional, simply connected hyperbolic space ${\bf H}^n$ of constant sectional curvature $-1$.}\\

As long as the isometry group of a Riemannian orbifold is
concerned, quite recently, A. V. Bagaev and N. I. Zhukova
\cite{BZ} showed the same result as Facts 1.1-2. They generalized
the idea of Kobayashi to their setting by using the orthonormal
frame bundle of a Riemannian orbifold. In this paper we consider a
special class of orbifolds --- manifolds with boundary. We firstly
observe that the dimension of the isometry group $I(M)$ of a
Riemannian manifold $M$ with boundary does not exceed
$\frac{1}{2}\,\dim\,M\,(\dim\,M-1)$. Then we classify such
Riemannian manifolds $M$ with boundary that the isometry groups
$I(M)$ attain the preceding maximal dimension. We divide the
lengthy classification list into three parts: Theorems 1.1-3, due
to that
their proofs will use different ideas. The notations in Fact 1.2 will be used in the following theorems.\\

\nd {\bf Theorem 1.1} {\it Let $M$ be an $n$-dimensional compact,
connected smooth Riemannian manifold with boundary and $n\geq 2$.
Suppose that the isometry group $I(M)$ is of dimension $n(n-1)/2$.
Then $M$ is diffeomorphic to either of the following four
manifolds\,{\rm :} the closed $n$-dimensional unit ball
$\overline{D^n}=\{x\in {\bf R}^n:|x|\leq 1\}$ in ${\bf R}^n$, the
two cylinder-like manifolds $S^{n-1}\times [0,\,1]$ and ${\bf
R}P^{n-1}\times [0,\,1]$, and the manifold ${\bf R}P^{n}\backslash
U$ constructed from ${\bf R}P^{n}$ with an $n$-dimensional open
disk $U\subset {\bf R}P^{n}$ removed, where the closure
$\overline{U}$ in ${\bf R}P^n$ is diffeomorphic to
$\overline{D^n}$. Furthermore, we can characterize the metric
tensor $g_M$ of $M$ as follows\,{\rm :}

\nd {\rm (1)}  If $M$ is diffeomorphic to $\overline{B^n}$, then
the metric $g_M$ of $M$ is rotionally symmetric with respect to a
unique interior point $O$ of $M$. That is, $g_M$ can be expressed
by $g_M=dt^2+\varphi^2(t)\, g_{S^{n-1}}$, where $g_{S^{n-1}}$ is
the standard metric on the unit sphere $S^{n-1}\subset {\bf R}^n$
and the function $\varphi:(0,\,R]\to (0,\,\infty)$ is smooth,
$\varphi(0)=0$, and
\[\varphi^{({\rm even})}(0)=0,\quad {\dot \varphi}(0)=1.\]

\nd {\rm (2)} If $M$ is diffeomorphic to $S^{n-1}\times [0,\,1]$,
then the metric $g_M$ can be expressed by
$dt^2+f^2(t)\,g_{S^{n-1}}$, where $T$ is a positive number and $f$
is a positive smooth function on $[-T/2,\,T/2]$.  The similar
statement holds for $M$ diffeomorphic to ${\bf R}P^{n-1}\times
[0,\,1]$.

\nd {\rm (3)} Suppose that $M$ is diffeomorphic to ${\bf
R}P^n\backslash U$. Then we can find a Riemannian manifold
$M'=S^n\times [-T/2,\, T/2]$ endowed with the metric
$dt^2+f^2(t)\,g_{S^{n-1}}$, where $f:[-T/2,\, T/2]\to
(0,\,\infty)$ is an even smooth function, and
 an involutive isometry $\beta$ of $M'$ defined by $\beta(x,\,t)=(-x,\,-t)$ such that $M$ is the quotient
 space of $M'$ by the group $\{1,\,\beta\}$.} \\

\nd {\bf Theorem 1.2} {\it Let $M$ be a noncompact connected
Riemannian manifold with boundary $\p M$ and of dimension $n\geq
2$ such that $\dim\, I(M)=n(n-1)/2$ and $\p M$ has at least one
compact component. Then $M$ is diffeomorphic to either
$S^{n-1}\times [0,\,1)$ or ${\bf R}P^{n-1}\times [0,\,1)$. In the
former case, the metric $g_M$ of $M$ can be expressed by
$dt^2+f^2(t)\,g_{S^{n-1}}$, where $f:[0,\,T)\to (0,\,\infty)$ is a
smooth function and $T$ is a positive number or $\infty$. Moreover
$M$ is complete if and only if $T=\infty$.
The similar statement holds for the latter case.}\\

\nd {\bf Theorem 1.3} {\it Let $M$ be a connected Riemannian
manifold with noncompact boundary $\p M$ and of dimension $n\geq
2$ such that $\dim\, I(M)=n(n-1)/2$. Denote by ${\bf H}^k$, $k\geq
2$, the $k$-dimensional complete simple connected Riemannian space
of constant sectional curvature $-1$. Then $M$ is diffeomorphic to
either ${\bf R}^{n-1}\times [0,\,1]$ or ${\bf R}^{n-1}\times
[0,\,1)$.
 Furthermore, we can characterize the metric tensor
$g_M$ of $M$ as follows\,{\rm :}

\nd {\rm (1)} If $M$ is diffeomorphic to ${\bf R}^{n-1}\times
[0,\,1]$, then there exists a positive number $T$ and a smooth
function $f:[0,\,T]\to (0,\,\infty)$ such that the metric tensor
$g_M$ on $M$ can be expressed by $g_M=dt^2+f^2(t)g_{{\bf
R}^{n-1}}$ or $g_M=dt^2+f^2(t)g_{{\bf H}^{n-1}}$ with $t\in
[0,\,T]$. Of course, we identify ${\bf H}^1$ with ${\bf R}^1$. The
metric $g_M$ in this case is always complete.

\nd {\rm (2)} If $M$ is diffeomorphic to ${\bf R}^{n-1}\times
[0,\,1)$, then there exists a number $T\in (0,\,\infty]$  and a
smooth function $f:[0,\,T)\to (0,\,\infty)$ such that the metric
tensor $g_M$ on $M$ can be expressed by $g_M=dt^2+f^2(t)g_{{\bf
R}^{n-1}}$ or $g_M=dt^2+f^2(t)g_{{\bf H}^{n-1}}$ with $t\in
[0,\,T)$. Moreover, the
metric $g_M$ is complete if and only if $T=\infty$.}\\

This paper is organized as follows. In section 2, we prove the
fact that the two definitions of isometry coincide on Riemannian
manifolds with boundary (see Proposition 2.1). It seems that
Bagaev-Zhukova \cite{BZ} did not mention this fact in their
setting of Riemannian orbifolds. The idea of making reduction to
the boundary in the proof of Proposition 2.1 will be used many
times afterwards. In this section we also show the above mentioned
observation  that the isometry group $I(M)$ of a Riemannian
manifold $M$ with boundary is a Lie transformation group of
dimension at most $\frac{1}{2}\,\dim\,M\,(\dim\,M-1)$ (see Theorem
2.1). Although our proof of the observation is based on the idea
of Proposition 2.1, to avoid the troublesome argument of point set
topology, we also use the result in Bagaev-Zhukova \cite{BZ} that
$I(M)$ has a Lie group structure. In sections 3 through 5, we use
the metric geometry and the theory of transformation group to
prove Theorems 1.1 through 1.3.

%%%%%%%%%%%%%%%%%%%%%%%%%%%%%%%% PRELIMINARIES %%%%%%%%%%%%%%%%%%%%%%%%%%%

\section{Some properties of isometry group}
\label{sec:def}

In the following sections we always let $M$ be an $n$-dimensional
connected, smooth Riemannian manifold with boundary and $n\geq 2$.
With the induced metric from $M$, the boundary $\p M$ of $M$ is an
$(n-1)$-dimensional Riemannian manifold without boundary. Note
that $\p M$ has at most countable connected components. Consider a
diffeomorphism $\phi$ of $M$ onto itself which preserves the
metric tensor. If $p$ is an interior point of $M$, then $\phi$
maps $p$ to another interior point, say $q$, and the differential
map $D\phi$ at $p$ induces an orthogonal transform from the
tangent space at $p$ to the one at $q$. If $u$ is a point on the
boundary $\p M$, then the tangent space $T_u M$ at $u$ should be
thought of as the upper half space
$$\{x=(x_1,x_2,\cdots,x_n)\in {\bf R}^n:x_n\geq 0\}$$
of the Euclidean space ${\bf R}^n$. That is,
\[T_u M=T_u (\p M)+\{\lambda{\bf n}_u:\lambda\geq 0\},\]
where ${\bf n}_u$ is the inner unit normal vector at $u$. Since
$\phi$ maps $u$ to another point $v$ on $\p M$, the differential
map $D\phi$ at $u$ maps ${\bf n}_u$ to ${\bf n}_v$, and maps
$T_u(\p M)$ orthogonally onto $T_v (\p M)$. Hence, in this sense,
we may also call that the differential map $D\phi$ at $u$ is an
orthogonal transform from $T_u M$ onto $T_v M$. Hence $\phi$
leaves the boundary $\p M$ invariant and induces an isometry of
$\p M$. Remember that by Myers-Steenrod \cite{MS} the two
definitions of isometry on $\p M$ are equivalent.

Let $d(\cdot,\,\cdot)$ be the distance function on $M$ induced by
the metric tensor of $M$ and $\psi$ a bijection on $M$ which
preserves $d(\cdot,\,\cdot)$. Since $\psi$ is a homeomorphism of
$M$ onto itself, its restrictions to the boundary $\p M$ is a
homeomorhism onto itself, so is the restriction to the interior of
$M$. In fact, we have a stronger property about $\psi$ in the
following

\begin{prop}
\label{prop:def} A distance-preserving bijection $\psi$ of $M$ is
a diffeomorphism which preserves the metric tensor of $M$. That
is, the two definitions of isometry of $M$ are equivalent.
\end{prop}
\nd \pf We firstly consider the property of $\psi$ near a point
$p$ of in the interior ${\rm Int}(M)$ of $M$. There exists an open
neighborhood $U\subset \im$ of $p$ such that the restriction
$\psi|_U$ of $\psi$ to $U$ is a distance-preserving map onto the
open neighborhood $V=\psi(U)\subset \im$ of $q=\psi(p)$. Since the
two definitions of isometry for Riemannian manifolds without
boundary are equivalent (see pp. 169-172 in Kobayashi-Nomizu
\cite{KN} for a proof), $\psi|_U:U\to V$ is a diffeomorphism
preserving the metric tensor. Hence, $\psi|_{\im}: \im\to \im$ is
also a diffeomorphism preserving the metric tensor. It suffices to
prove that for each point $p\in \p M$, $\psi$ is smooth near $p$
and $D\psi$ at $p$ is an orthogonal transform from $T_p M$ onto
$T_q M$, where $q=\psi(p)\in \p M$. We divide the proof into two
steps.

{\it Step 1}\quad Remember that $\bn_p$ denotes the  inner unit
normal vector at $p$. Choose $\d>0$ so small that the  geodesic
$\gamma(p,\,t):=\exp_p\,(t{\bf n}_p)$, $t\in [0,\,\d]$, satisfies
\begin{equation}
\label{equ:21} d\bigl(\g(p,\,t),\, \p
M)=d\bigl(\g(p,\,t),\,\g(p,\,0)\bigr)=t.
\end{equation}
Since the geodesic emanating from each point with this property is
unique, the image $\psi\circ\g$ of $\g$ under the isometry $\psi$
is also a geodesic perpendicular to the boundary $\p M$ at the
initial point $q$. Actually, we will see later that the
differential $D\psi$ at $p$ maps $\bn_p$ to $\bn_q$. In $\p M$, we
choose a small open neighborhood $V\subset \p M$ of $p$ such that
for each point $p'\in V$ the geodesic $\gamma_(p',\,t)$, $t\in
[0,\,\d/2]$, satisfies (\ref{equ:21}). Then the map
$\g(\cdot,\,t):p'\mapsto \g(p',\,t)$ gives a diffeomorphism of $V$
onto a hypersurface $V_t$ in $\im$ for each $t\in (0,\,\d/2]$. We
can define the similar map from the neighborhood $\psi(V)\subset
\p M$ of $q$ and denote the map also by $\g(\cdot,\,t)$. We
observe that
\[\psi\circ \g(\cdot,\,t)=\g(\cdot,\,t)\circ\psi\]
holds for each $p\in V$ and each $t\in [0,\,\delta/2]$ so that
$\g(\psi(V),\,t)=\psi(V_t)$. Hence the map $\psi|_V$  can be
thought of as the composition of three diffeomorphisms,
\[\psi|_V(\cdot)= \bigl(\g(\cdot,\,t)\bigr)^{-1}\circ\psi\circ\g(\cdot,\,t).\]
So $\psi|_V$ is a diffeomorphism of $V$ onto $\psi(V)$. Since
$\psi|_{V_t}$ is an isometry of $V_t$ onto $\psi(V_t)$ for each
$0<t\leq \d/2$, letting $t\to 0$, we can see that $\psi|_V$ is a
one-to-one distance-preserving map of $V$ onto $\psi(V)$. Since
the two definitions of isometry for Riemannian manifolds without
boundary are equivalent,
 $\psi|_V$ is a diffeomorphism onto $\psi(V)$ preserving the metric tensor.

{\it Step 2}\quad By Step 1 and its preceding argument, we can
take a small open neighborhood $U\subset M$ of $p$ such that the
exponential map $\exp_p$ at $p$ is a diffeomorphism from some
neighborhood ${\tilde U}$ of $0$ in $T_p M$ onto $U$. Remember
that the partial derivative of $\psi$ exists at the direction of
the inner unit normal ${\bf n}_p$, and equals ${\bf n}_{q}$. That
is, the differential map $D\psi$ at each $p$ in $U\cap \partial M$
is an orthogonal transform from $T_p M$ to $T_{\psi(p)}M$. Since
$\psi$ is a homeomorphism, we may assume that $\psi(U)$ is
contained in a normal coordinate neighborhood of $q$. Since both
$\exp_p$ and $\exp_q$ are local diffeomorphisms, the equality
\[\psi\circ\exp_p=\exp_q\circ (D\psi(p))\]
in ${\tilde U}$, which gives that $\psi|_U$ is diffeomorphism
onto $\psi(U)$. \hf \\

The above proof essentially follows the idea in pp. 169-172 of
Kobayashi-Nomizu \cite{KN}. We repeat it here because its idea and
notations will be used many times later. The following lemma is
elementary and useful, but its proof is omitted.

\begin{lem}
\label{lem:id} Let $\phi$ be an isometry of $M$. If $\phi$ has a
fixed point $p\in M$ and the differential $D\phi$ at $p$ is the
identity map of $T_p M$, then $\phi$ is the identity map of $M$.
\end{lem}

\begin{prop}
\label{prop:closed} \nd {\rm (1)} Let $\phi$ be an isometry of
$M$. Then the restriction $\phi|_{\p M}$ to the boundary $\p M$ is
an isometry of $\p M$. Moreover, if $\phi$ leaves each point of a
component $B$ of $\p M$ fixed, then $\phi$ is the identity map of
$M$.

\nd{\rm (2)} Let $\phi$ be an element in the identity component
$I^0(M)$ of the isometry group of $M$. The restriction $\phi|_B$
to a connected component $B$ of $\p M$ is an element in the
identity component $I^0(B)$ of the isometry group of $B$. This map
$$\iota:I^0(M)\to I^0(B),\quad \phi\mapsto \phi|_B$$
gives a continuous monomorphism with the image closed  in
$I^0(B)$. That is, $\iota$ is a regular embedding of the Lie group
$I^0(M)$ into the Lie group $I^0(B)$.
\end{prop}

\nd\pf

\nd (1) The first statement have been shown in the proof of
Proposition \ref{prop:def}. If each point $p$ of a component $B$
of $\p M$ is fixed by $\phi$, then the differential $D\phi$ at $p$
is the identity map of $T_p M$. By Lemma \ref{lem:id} $\phi$ is
the identity map of $M$.

\nd (2) Given a point $p$ in a connected component $B$ of
$\partial M$, we claim that each element $\phi$ in $I^0(M)$ maps
$p$ to a point in $B$. Otherwise, we assume that $\phi(p)$ lies in
another component $B'$ distinct from $B$. Choosing a path
$\{\phi_t\}$ in $I^0(M)$ with $\phi_0=id_M$ and $\phi_1=\phi$, we
obtain a path $\{\phi_t(p)\}$ connecting $p$ and $\phi(p)$. Since
each diffeomorphism of $M$ maps $\partial M$ onto $\partial M$, we
find that the path $\{\phi_t\}$ lies on $\partial M$ and reach a
contradiction. By the proof of Proposition
 \ref{prop:def}, we know that $\phi|_B$ is an isometry of $B$, which is actually an element of $I^0(B)$.
Moreover, if $\phi|_B$ is the identity map of $B$, then by (1)
$\phi$ is the identity map of $M$. Thus, the map $\iota:I^0(M)\to
I^0(B),\ \phi\mapsto \phi|_B$ gives a continuous monomorphism of
$I^0(M)$ into $I^0(B)$.

Finally we need to show that the image of $I^0(M)$ under $\iota$
is closed in $I^0(B)$ with respect to the compact-open topology.
We divide the proof into two steps.

{\it Step 1} We show that for a sequence $\{\phi^n\}$ of
isometries in $I^0(M)$ such that
\[\phi^n|_B\to id_B\quad {\rm in}\quad I^0(B),\]
there holds $\phi^n\to id_M$ in $I^0(M)$. By Fact 1.1,  $I^0(B)$
has the structure of a Lie group, whose topology from its smooth
structure coincides with the compact-open topology. We may assume
that $\phi^n|_B\to id_B$ in the $C^1$ topology of $I^0(B)$. That
is, $\phi^n$ converges to the identity map $id_K$ in the sense of
the $C^1$ norm in each compact neighborhood $K$ in the topological
space $B$. Then, by using the normal coordinate charts with
respect to each point in $K$ (see the last three lines of the
proof of Proposition 2.1), we find that there exists a compact
neighborhood ${\cal K}$ of $M$ such that $K\subset {\cal K}$ and
$\phi^n(p)$ converges uniformly to $p$ for each point $p\in{\cal
K}$.

{\it Step 2} We show that if a sequence $\{\phi^n\}$ of isometries
in $I^0(M)$ satisfies $\phi^n|_B\to \phi$ in $I^0(B)$, then there
exists $\psi\in I^0(M)$ such that $\psi|_B=\phi$. Since, by
Myers-Steenrod \cite{MS} and Bagaev-Zhukova \cite{BZ}, both
$I^0(B)$ and $I^0(B)$ have structures of Lie groups, they can be
endowed with a Riemannian metric. By the Cauchy criterion,
$\phi^n|_B$ converges in $I^0(B)$ $\phi^n|_B(\phi^m|_B)^{-1}$
converges to $id_B$ as $m,n\to \infty$. Step 1 tells us that
$\phi^n(\phi^m)^{-1}$ converges to $id_M$ as $m,n\to \infty$. That
is,
$\phi^n$ converges to some $\psi\in I^0(M)$ such that $\psi|_B=\phi$.  \hf \\

As an immediate corollary of Proposition 2.2 and Fact 1.3, we obtain\\

\nd {\bf Theorem 2.1.} {\it The isometry group $I(M^n)$ has a
structure of Lie group of dimension at most
$n(n-1)/2$.}\\

%%%%%%%%%%%%%%%%%%%%%%%%%%%%% PROOF OF THEOREM CUP  %%%%%%%%%%%%%%%%%%%%%%%%%%

\section{Proof of Theorem 1.1}
Let $M$ be a Riemannian manifold satisfying the condition of
Theorem 1.1 in this section. By Proposition 2.2 and Fact 1.2, the
isometry group $I(B)$ of each component $B$ of $\p M$ attains the
maximal dimension $n(n-1)/2$, so $B$ is isometric to either
$S^{n-1}$ or ${\bf R}P^{n-1}$ with constant sectional curvature
$1$. If $n=2$, the $\partial M$ consists of circles. But our
argument later also goes through in this case.

Suppose that there exists a component $B$ isometric to the sphere
$S^{n-1}$. By the proof (see pp. 46-47 in \cite{Ko}) of Fact 1.2,
$G:=I^0(M)$ is isomorphic to $SO(n)$ and its action on $B$ is just
the linear action of $SO(n)$ on $S^{n-1}$. We may identify $G$
with $SO(n)$, $B$ with $S^{n-1}$ up to a scaling of metric.
Remember that $G$ acts transitively on $S^{n-1}$ and the isotropy
group $G_x$ at each point $x$ of $S^{n-1}$ is isomorphic to
$H:=SO(n-1)$. Here we use the notation in the proof of Proposition
2.1. Choose a positive number $\d>0$ such that the map
$\g(\cdot,\,t)$ is a diffeomorphism of $B$ onto a hypersurface
$B_t$ in $\im$ for each $t\in (0,\,\d]$. Since the $G$-action
interchanges with $\g(\cdot,\,t)$, $G$ leaves each $B_t$
invariant. We claim that the $G$-action on $M$ has the principal
orbit whose type is $SO(n)/SO(n-1)$.  Remember that  the union of
the principal orbits of type different from $G/H$ forms an open
and dense subset of $M$. For the detail of this, see Theorem 4.27
in pp. 216-220 of Kawakubo \cite{Ka}, where only manifolds without
boundary are considered. But Theorem 4.27 in Kawakubo \cite{Ka}
also holds for our case by virtue of the map $\g(\cdot,\,t)$. The
claim follows from that the union $M(H)$ of orbits with type $G/H$
contains the open subset $\cup_{t\in [0,\,\d]}\, B_t$ in $M$.
Therefore, every component of $\p M$ is isometric to $S^{n-1}$. By
Theorem 4.19 (pp. 202-203) and  Theorem 4.27 in Kawakubo
\cite{Ka}, the orbit space $M(H)/G$ is a connected smooth
1-manifold, whose boundary coincides with the orbit space $\p
M/G$. Hence, $\p M$ has at most two components. We also have the
similar argument when $\p M$ has a component isometric
to ${\bf R}P^{n-1}$. Summing up, we have proved \\

\nd {\bf Lemma 3.1} {\it  $\p M$ has at most two components, each
of which is isometric to either $S^{n-1}$ or ${\bf R}P^{n-1}$.
Moreover, if $\p M$ has two components, then the two components
are isometric up to a scaling
of metric.}  \\

\nd{\bf Lemma 3.2} {\it Let $\p M$ have two components. Then $M$
is diffeomorphic to either $S^{n-1}\times [0,\,1]$ or ${\bf
R}P^{n-1}\times [0,\,1]$. If $M$ is diffeomorhic to $S^{n-1}\times
[0,\,1]$, then there exists a positive number $T$ and a positive
smooth function $f(t)$ on $[0,\,T]$ such that the Riemannian
metric $g_M$ on $M$ can be expressed by
\[g_M=dt^2+f^2(t)\,g_{S^{n-1}}.\]
The similar statement holds for $M$ diffeomorphic to ${\bf R}P^{n-1}\times [0,\,1]$. }\\

\nd \pf We only prove the case that both the two component of $\p
M$ are isometric to $S^{n-1}$.
 By the proof of Lemma 3.1, the orbit space $M(H)/G$ is a
closed interval. On the other hand, the orbit space $M(H)/G$ is a
dense subset in the total orbit space $M/G$. Therefore, $M=M(H)$,
i.e. all the orbits of $G$-action on $M$ is of principal type. So
$M$ is a smooth fiber bundle with fiber $S^{n-1}$ over a compact
interval. Actually $M$ is diffeomorphic to the product of
$S^{n-1}$ and a compact interval.

We use the notation in the proof of Proposition \ref{prop:def} in
what follows. Let $B$ and $B'$ be the components of $\p M$. The
two spheres $B$ and $B'$ may have different size. Choose a point
$p\in B$. Then we claim that the geodesic $\g(p,\,t)$ with initial
velocity $\bn_p$ terminates at some point $p'\in B'$ and with the
ending verlocity $\bn_{p'}$ at a positive time determined later.
Indeed, choosing $p'\in \p M$ such that
$d(p,\,B')=d(p,\,p')=:T(p)$, we can find a geodesic between $p$
and $p'$ whose length is $T(p)$. It is clear that this geodesic is
perpendicular to the boundary $\p M$ at $p'$. It also is
perpendicular to $\p M$ at $p$. Otherwise, we can find another
point $q$ on $B$ and a path $\ell$ connecting $q$ and $p'$ of
length less than $T(p)$. Choosing an element $\al\in G$ mapping
$q$ to $p$, we obtain a path $\al(\ell)$ connecting $p$ and
$\al(p')\in B'$ of length less than $T(p)$. Contradiction!
Therefore, this geodesic is exactly the one $\{\g(p,\,t):t\in
[0,\,T(p)]\}$ in the claim. We claim again that $T(p)$ equals the
distance $T$ bewteen $B$ and $B'$. Actually, there exists $p_0$
such that the geodesic $\g(p_0,\,[0,\,T(p_0)])$ satisfies
\[T(p_0)=d\bigl(p_0,\,\g(p_0,\,T(p_0)\bigr)=d(B,\,B').\]
But $G$ acts transitively on the set of geodesics
$$\left\{\g\bigl(p,\,[0,\,T(p)]\bigr):p\in B\right\}$$
of geodesics connecting $B$ and $B'$, which implies that these
geodesics have the same length $T$. So each point $p$ in $B$ has
the same distance $T$ to $B'$.

We claim that the map $\g(\cdot,\,t):B\to M$ is a diffeomorphism
of $B$ onto the hypersurface $B_t$ for each $t\in [0,\,T]$.
Similarly as the preceding paragraph, this map is sujective.  It
is also injective. Otherwise, there exists $0<T'<T$ and two
distinct points $p_1$ and $p_2$ in $B$ such that $\g(p_1,\,
T')=\g(p_2,\,T')=:q$. Then we get a curve $\g(p_1,\,[0,\,T'])\cup
\g(p_2,\,[T',\, T])$ irregular at point $q$ with length $T$
connecting $B$ and $B'$. Contradiction. Since this map is
equivariant with respect to the $G$ actions on $B$ and $B_t$,
every point of $B_t$ is a regular point of the map by the Sard
theorem. Combining above, we know that the map $\g(\cdot,\,t):B\to
B_t$ is one-to-one, onto and its differential does not degenerate.
That is, it is a diffeomorphism.

Then we show that any two geodesics in the set
$\left\{\g\bigl(p,\,[0,\,T]\bigr):p\in B\right\}$ are equal if
they intersect. Suppose that there exists distinct points
$p_1,\,p_2\in B$ and positive times $t_1,\,t_2$ in $(0,\,T)$ such
that $\g(p_1,\,t_1)=\g(p_2,\,t_2)$. Then by the argument in the
precedent paragraph we know $t_1\not=t_2$, say $t_1<t_2$. Then the
piecewise smooth curve $\g(p_1,\,[0,\,t_1])\cup \g(p_2,\,[t_2,\,
T])$ connecting $B$ and $B'$ has length $t_1+(T-t_2)<T$.
Contradiction.

The statements in the last two paragraphs shows that the map
\[\g(\cdot,\,\cdot):[0,\,T]\times B\to M,\quad (p,\,t)\mapsto \g(p,\,t)\]
is a diffeomorphism. Moreover, since $G$ leaves each hypersurface
$B_t$ invariant and acts isometrically and effectively on it,
$B_t$ is a sphere with constant sectional curvature. Therefore,
the Riemannian metric $g_M$ of $M$ can be written by
$g_M=dt^2+f^2(t)\, g_{S^{n-1}}$ for some positive smooth function
$f(t)$ on $[0,\,T]$.
 \hf\\

\nd{\bf Lemma 3.3} {\it If $\p M$ is connected, then $\p M$ must be isometric to the unit sphere $S^{n-1}$.}\\

\nd \pf Suppose that $\p M$ is isometric to the real projective
space ${\bf R}P^{n-1}$ with the constant sectional curvature and
$n\geq 3$. Then $n$ should be even since ${\bf R}P^{\rm even}$
does not bound by the unoriented cobordism theory (see pp. 52-53
in Milnor-Stasheff \cite{MS}). Denote by $G$ the identity
component of the isometry group of $M$. Since $\dim\, G=n(n-1)/2$,
$G$ is isomorphic to $SO(n)/\{\pm 1\}$ and its action on the
boundary ${\bf R}P^{n-1}$ is induced by the linear action of
$SO(n)$ on $S^{n-1}$. Moreover, the isotropy subgroup $H:=G_p$ at
each point $p$ on $\p M$ is isomorphic to $SO(n-1)$. Following the
proof of Lemma 3.1, the $G$ action on $M$ has principal orbits of
type $G/H={\bf R}P^{n-1}$. Denote by $M(H)$ the union of principal
orbits. Then the orbit space $M(H)/G$ is diffeomorphic to the
interval $[0,\,1)$.  Since $M(H)/G$ is dense in the total orbit
space $M/G$, for the $G$ action on $M$, there exists only one
orbit $G/J$ other than the principal ones. We call the orbit $G/J$
the exceptional orbit. It is this orbit $G/J$ that corresponds to
the endpoint 1 of the orbit space $M/G$. Since $H$ can be thought
of as a proper subgroup of $J$ in the sense of conjugacy, the Lie
algebra ${\frak J}$ of $J$ contains a subalgebra isomorphic to
${\frak so}(n-1)$. Simple computation shows that ${\frak J}$ is
isomorphic to either ${\frak so}(n)$ or ${\frak so}(n-1)$.

{\it Case 1} If ${\frak J}$ is ${\frak so}(n)$, then $J$ equals
$G$. This means that topologically $M$ is the cone of ${\bf
R}P^{n-1}$, which does not have the structure of a manifold.
Actually, the cone of ${\bf R}P^{n-1}$ is homeomorphic to the
orbifold $\overline{D^n}/\{\pm 1\}$ the boundary ${\bf R}P^{n-1}$.
Contradiction.

{\it Case 2} If ${\frak J}$ is isomorphic to ${\frak so}(n-1)$,
then the exception orbit $G/J$ is a submanifold of codimension 1
in $M$. So the $G$ action on $G/J$ is also effective. Since $G$
attains the largest-possible dimension $n(n-1)/2$,  $G/J$ is
isometric to ${\bf R}P^{n-1}$ by Fact 1.3, which implies that
$G/J$ is also a principal orbit. Contradiction.
\hf\\

\nd {\bf Lemma 3.4} {\it If $\p M$ is isometric to the sphere
$S^{n-1}$, then $M$ is homeomorphic to either $\overline{D^n}$ or
${\bf R}P^n\backslash U$, where $U$ is an $n$-dimensional open
disk such that its closure $\overline{U}$ in ${\bf R}P^n$ is
homeomorphic to the closed unit ball $\overline{B^n}\subset
{\bf R}^n$.}\\

\nd \pf We use the notation in the proof of Lemmas 3.1-3. Remember
that $G=SO(n)$ and $H=SO(n-1)$. By the similar argument in the
proof of Lemma 3.3, there exists only one orbit $G/J$ (called the
exceptional orbit) other than the principal orbits among all the
orbits of the $G$ action on $M$. Also by the argument in the proof
of Lemma 3.3, we know that either $J$ is $SO(n)$ or $J$ contains a
subgroup of finite index and isomorphic to $SO(n-1)$.

{\it Case 1} If $J$ is $SO(n)$, then $M$ is homeomorphic to the
closed unit disk.

{\it Case 2} If $SO(n-1)$ is a subgroup of $J$ of finite index,
then the exceptional orbit $G/J$ is a submanifold of codimension 1
in $M$. Hence $G$ acts effectively  on $G/J$. Since $\p M$ is
connected, $G/J$ has to be isometric to ${\bf R}P^{n-1}$ by Fact
1.3. Therefore, $M$ is homeomorphic to the mapping cone
$S^{n-1}\times [0,\,1]/\sim$, where the equivalent relation $\sim$
on $S^{n-1}\times [0,\,1]$ is defined by $(x,\,1)\sim (-x,\,1)$
(see p. 13 for the concept of mapping cone in p. 13 Hatcher
\cite{Ha}). Clearly $M$ is also
homeomorphic to the punctured real prjective space ${\bf R}P^n\backslash U$.\hf \\

\nd{\bf Lemma 3.5} {\it If $M$ is homeomorphic to the
$n$-dimensional closed disk $\overline{D^n}$, then there exists a
point $O$ in the interior of $M$ and a positive number $R>0$ such
that the exponential map $\exp_O$ at $O$ is a diffeomorphism of
the closed ball centered at the origin $0\in T_O M$ and of radius
$R$ in $T_O M$ onto $M$. Moreover,  the Riemannian metric $g_M$ of
$M$ is rotationally symmetric with respect to $O$ so that it can
be expressed by
\[g_M=dt^2+\varphi^2(t)\, g_{S^{n-1}},\]
where the function $\varphi:(0,\,R]\to (0,\,\infty)$ is smooth,
$\varphi(0)=0$, and
\[\varphi^{({\rm even})}(0)=0,\quad {\dot \varphi}(0)=1.\]
}

\nd \pf We use the notation in the proof of Theorem 2.1 and Lemma
3.4. Denote by $B$ the boundary $\p M$. We know that $B$ is a
sphere of constant sectional curvature. By the proof of Lemma 3.4,
there exists a unique fixed point $O\in {\rm Int}(M)$ of the $G$
action on $M$. The $G$ action on the tangent space $T_O M$ at $O$
gives an isomorphism of $G$ onto the special orthogonal
transformation group $SO(T_O M)$ of the Euclidean space $T_O M$.
We denote by $B(0,\, r)$ the set of tangent vectors of length $<r$
at $O$, and by $S(0,\,r)$ the set of tangent vectors of length $r$
at $O$. We also denote by $S(r)$ the set of points in $M$ with
distance $r$ to $O$, by $\exp$ the exponential map at $O$.

Choose a point $p\in B$ such that $d(O,\,p)=d(O,\,B)=:R$. Then
there exists a unit tangent vector $V_0$ at $O$ such that
$p=\exp\, (RV_0)$, and the geodesic $\exp\,(tV_0)$, $t\in
[0,\,R]$, is perpendicular to $B$ at $p$. Choose an arbitrary unit
tangent vector $V$ at $O$. We claim that the geodesic $\exp\,
(tV)$, $t\geq 0$, meets the bounadry $B$ perpendicularly at time
$R$. Actually, choosing an isometry $\alpha\in G$ such that the
differential $d\alpha$ at $O$ maps $V_0$ to $V$, we find that the
geodesic $\exp_O\, (tV)$, $t\in [0,\,R]$, is the image of the one
$\exp\, (tV_0)$, $t\in [0,\,R]$, under the isometry $\alpha$.
Since $G$ acts transitively on $B$, for each point $p\in B$, there
exists $V\in S(0,\,1)$ such that $\exp\,(RV)=p$. By the uniqueness
of the geodesic perpendicular to $B$ at a given point, we find
that if $\exp\,(RV)=\exp\,(RW)$ for any two vectors $V,\,W\in
S(0,\,1)$, then $V=W$. Hence $\exp:\, S(0,\,R)\to S(R)$ is a
smooth bijection and there is no cut point of $O$ in the interior
of $M$. So, $\exp$ gives a diffeomorphism of $B(0,\,R)$ onto ${\rm
Int}(M)$. To prove this diffeomorphism can extend to the boundary,
by the Gauss lemma, we need only to show that the restriction of
$\exp$ to $S(0,\,R)$ is a diffeomorphism onto $S(R)=B$.

Remember that $G$ acts on both $S(0,\,R)$ and $S(R)$. Moreover,
the exponential map $\exp:\, S(0,\,R)\to S(R)$ is an equivariant
smooth bijection map with respect to the $G$ actions. By the Sard
theorem, there exists a regular value of $\exp|_{S(0,\,R)}$. On
the other hand, by the equivariant property, all points of $S(R)$
are regular values. That is, the map $\exp:\, S(0,\,R)\to S(R)$ is
a diffeomorphism.

The similar statement holds  for each $(n-1)$ dimensional sphere
$S(r)$, $0<r\leq R$, in $M$. So the metric $g_M$ of $M$ has
rotational symmetry with respect to $O$. The expression
of $g_M$ follows from the argument in pp. 12-13 in Petersen \cite{Pe}.\hf \\

\nd{\bf Remark 3.1} We can list closed geodesic balls with
suitable radii in the three spaces ${\bf R}^n$, $S^n$ and $H^n$ as
concrete examples of the manifold $M$ in Lemma 3.5.
Simultaneously, the geodesic annuli of these three spaces form
examples of the manifold in Lemma 3.4. The manifold $M$ in our
consideration need not have constant sectional curvature, whose
curvature can be computed explicitly in terms of the function $f$
(see pp. 65-68 in Petersen \cite{Pe}). Because of the large
symmetry on them, this class of manifolds, including geodesic
balls in ${\bf R}^n$ , $S^n$
and $H^n$, may be thought of as the simplest class of compact Riemannian manifolds with boundary.  \\

\nd {\bf Lemma 3.6} {\it We use the notation in the proof of Lemma
3.4. Suppose that $M$ is homeomorphic to ${\bf R}P^n\backslash U$.
Then we can find a Riemannian manifold $M'=S^n\times [-T/2,\,
T/2]$ endowed with the metric $dt^2+f^2(t)\,g_{S^{n-1}}$, where
$f:[-T/2,\, T/2]\to (0,\,\infty)$ is an even smooth function, and
an involutive isometry $\beta$ of $M'$ defined by
$\beta(x,\,t)=(-x,\,-t)$ such that $M$ is the quotient space of
$M'$ by the group $\{1,\,\beta\}$. Here $-x$ means the antipodal
point of $x$ in $S^{n-1}$.
Of course, $M$ is diffeomorphic to ${\bf R}P^n\backslash U$.}\\

\nd \pf First of all, let us forget the Riemannian metric on $M$.
Consider a topological model of $M$ --- the mapping cone
$S^{n-1}\times [0,\,1]/\sim$. Remember that the equivalent
relation $\sim$ means $(x,\,1)\sim (-x,\,1)$, where $x\mapsto -x$
is the deck transformation of the 2-fold covering $S^{n-1}\to {\bf
R}P^{n-1}$. Then  $M$ is the quotient of  $M':=S^{n-1}\times
[0,\,2]$ by the group generated by the involution $\beta$ of $M$
given by
\[\beta(x,\,t)=(-x,\,2-t).\]

Then we endow $M'$ with the induced Riemannian metric from $M$.
Since each isometry of $M$ can be lifted to two isometries of
$M'$, $M'$ also satisfies the condition of Theorem 1.3. By Lemma
3.2, there exists a positive number $T$ and a smooth function
$f:[-T/2,\,T/2]\to (0,\,\infty)$ such that $M'$ is diffeomorphic
to $S^{n-1}\times [-T/2,\,T/2]$ and the metric $g_{M'}$ is given
by
\[g_{M'}=dt^2+f^2(t)\,g_{S^{n-1}},\]
where $g_{S^{n-1}}$ is the standard metric on the unit sphere
$S^{n-1}$. On the other hand, since the deck transform
$\beta:M'\to M'$, $(x,\,t)=(-x,\,-t)$, is an isometry of $M'$,  we
can see that
$-x$ is actually  the antipodal point of $x\in S^{n-1}$ and $f(t)$ is an even function.\hf \\

We finally complete the proof of Theorem 1.1 by combining all the
lemmas in this section. \hf

%%%%%%%%%%%%%%%%%%%%%%%%%%%%%%%%%%%%%  Section 4

\section{Proof of Theorem 1.2}

Let $M$ be a Riemannian manifold satisfying the condition of
Theorem 1.2 in this section. Let $B$ be a compact component of $\p
M$. By Proposition \ref{prop:closed} and Fact 1.3, the isometry
group $I(B)$ of $B$ has the largest-possible dimension $n(n-1)/2$,
so $B$ is isometric to either $S^{n-1}$ or ${\bf R}P^{n-1}$ with
constant sectional curvature. Suppose the former case holds. Then
$G:=I^0(M)$ is isomorphic to $SO(n)$ and the isotropy subgroup
$G_x$ at each point $x\in B$ is isomorphic to $H:=SO(n-1)$. By the
same argument in the proof of Lemma 3.1, the $G$ action on $M$ has
principal orbit type $G/H$ and $\p M$ has at most two components.
We claim that $\p M$ is connected, i.e. $\p M=B$. Actually, since
the orbit space $M(H)/G$ is a connected smooth 1-manifold, whose
boundary coincides with the orbit space $\p M/G$, $M(H)/G$ is a
compact interval if $\p M$ have two components. Since the orbit
space $M(H)/G$ is dense in the total orbit space $M/G$, we have
$M(H)=M$ and $M$ is diffeomorphic to $S^{n-1}\times [0,\,1]$.
Contradict the noncompactness of $M$. By the same reason, we also
have $M(H)=M$. We have the similar argument for the latter case.
Summing up, we obtain\\

\nd {\bf Lemma 4.1} {\it $M$ is diffeomorphic to either
$S^{n-1}\times [0,\,1)$ or
${\bf R}P^{n-1}\times [0,\,1)$.}\\

\nd {\bf Lemma 4.2}

\nd (1) {\it Let $M$ be complete. If $M$ is diffeomorphic to
$S^{n-1}\times [0,\,1)$, then the metric $g_M$ of $M$ can be
expressed by
\[g_M=dt^2+f^2(t)\,g_{S^{n-1}},\]
where $f:[0,\,\infty)\to (0,\,\infty)$ is a  smooth function. The
similar statement holds for $M$ diffeomorphic to ${\bf
R}P^{n-1}\times [0,\,\infty)$.}

\nd (2) {\it Let $M$ be noncomplete. If $M$ is diffeomorphic to
$S^{n-1}\times [0,\,1)$, then there exists a positive number $T$
the metric $g_M$ of $M$ can be expressed by
\[g_M=dt^2+f^2(t)\,g_{S^{n-1}},\]
where $f:[0,\,T)\to (0,\,\infty)$ is a smooth function.
The similar statement holds for $M$ diffeomorphic to ${\bf R}P^{n-1}\times [0,\,1)$.} \\

\nd \pf We only prove the case that $M$ is diffeomorphic to
$S^{n-1}\times [0,\, 1)$. We use the notation in the proof of
Proposition \ref{prop:def} and Lemma 4.1. Choose an arbitrary
point $q$ in the interior of $M$ such that the distance of $q$ to
the boundary $B$ is $D$. Cutting $M$ along the orbit through $q$
of the $G$ action, we obtain a compact part $M_1$ diffeomorphic to
$S^{n-1}\times [0,\,1]$ and a noncompact part diffeomorphic to
$M$, on both of which $G$ acts isometrically. By the proof of
Lemma 3.2, the map $\g(\cdot,\,\cdot):B\times [0,\,D]\to M_1$ is a
diffeomorphism. Since $q$ is arbitrary, there exists $T\in
(0,\,\infty]$ such that the map $\g(\cdot,\,\cdot):B\times
[0,\,T)\to M$ is a diffeomorphism. Moreover, $T$ is $\infty$ if and only if $M$ is complete.\hf\\

%%%%%%%%%%%%%%%%%%%%% section 5
\section{Proof of Theorem 1.3}

Since the compact transformation group theory cannot be applied
directly to isometry groups of Riemannian manifolds with
noncompact boundary, we need new ideas to classify Riemannian
manifolds with noncompact boundary  whose isometry groups attain
the maximal dimension.

Denote by $G_k$ and ${\cal G}_k$, the identity components of the
isometry groups of ${\bf R}^k$ and ${\bf H}^k$, respectively.
Remember that ${\cal G}_k$ is the identity component of $O(1,\,k)$
and semisimple for each $k\geq 2$. However, $G_k$ is the
semidirect product of $SO(k)$ and ${\bf R}^k$, and it is not
semisimple for each $k\geq 1$ (see p. 5 and p. 77 in Petersen
\cite{Pe}). Let $M$ be a Riemannian manifold satisfying the
assumption of Theorem 1.5 through this section. By Proposition
\ref{prop:closed}  and Fact 1.2, every component of $\p M$ with
the induced Riemannian metric from $M$ is isometric to either
${\bf R}^{n-1}$ or the $(n-1)$-dimensional complete and simply
connected Riemannian manifold $H^{n-1}(c)$ of constant sectional
curvature $c<0$. Note that all $H^{n-1}(c)$'s, $c<0$, have the
same isometry group isomorphic to the semidirect product of ${\cal
G}_{n-1}$ and ${\bf Z}_2$. Suppose that a component of $\p M$ is
isometric to ${\bf R}^{n-1}$. Then $I^0(M)$ is isomorphic to
$G_{n-1}$, which acts effectively and isometrically on each
component of $\p M$. Hence, we find that each component of $\p M$
should be isometric to ${\bf R}^{n-1}$. The similar argument goes
through if $\p M$
has a component isometric to $H^{n-1}(c)$ for some $c<0$. \\

\nd {\bf Lemma 5.1.} {\it Each component of $\p M$ is isometric to
either ${\bf R}^{n-1}$ and $H^{n-1}(c)$ for
some $c<0$. Moreover, the components of $\p M$ are mutually isometric up to a scaling of metric.}  \\

\nd {\bf Lemma 5.2.} {\it We use the notation of Proposition
\ref{prop:def}. Let $B$ be a component of the boundary $\p M$ and
$p$ an arbitrary point of $B$. Let $I$ be the maximal existence
interval of the geodesic $\g(p,\,t)=\exp_{p} (t{\bf n}_p)$
perpendicular to $B$ at the initial point $p$. Then the map
\[\g: B\times I\to M,\quad (q,\,t)\mapsto \g(q,\,t)\]
is well defined and gives a diffeomorphism of $B\times I$ onto
$M$. Consequently, if $I$ is a compact interval, then $M$ is
diffeomorphic to $B\times [0,\,1]$; if $I$ is an interval
open at the right endpoint, then $M$ is diffeomorphic to $B\times [0,\,1)$.}\\

\nd\pf Denote by $G$ the identity component of $I(M)$. Then
$G=G_{n-1}$ if $B$ is isometric to ${\bf R}^{n-1}$, $G={\cal
G}_{n-1}$ if $B$ is isometric to $H^{n-1}(c)$ for some $c<0$. We
only prove the case where $B$ is isometic to ${\bf R}^{n-1}$. The
other case can be proved similarly.

Since $G$ acts transitively on $B$, we can see that, for every
point $q\in B$, the geodesic $\g(q,\,t)$ perpendicular to $B$ at
the initial point $q$ also has the maximal existence interval $I$.
We claim that for any two distinct points $p,\,q\in B$, the two
geodesics $\{\g(p,\,t):t\in I\}$ and $\{\g(q,\,t):t\in I\}$ does
not intersect at a point $x$ such that $x=\g(p,\,s)=\g(q,\,s)$ for
some $s\in I$ and $d(x,\,p)=d(x,\,q)=s$. Otherwise, we take an
element $\al\in G$ mapping $p$ to $q$ and the sequence of points
$p_n=\al^n(p)$ goes to the infinity of $B$ as $n\to\infty$, which
implies that the subgroup $\Gamma$ generated by $\al$ is not
precompact in $G$.  On the other hand, since
$\g\bigl(\al(p),\,s\bigr)=\al\bigl(\g(p,\,s)\bigr)$, $x$ is the
fixed point of $\al$. So the isotropy subgroup $G_x$ of the $G$
action at $x$ contains a non-precompact subgroup $\Gamma$.
Contradict the compactness of the isotropy group $G_x$.

We claim that the subset $B_t=\{\g(p,\,t):p\in B\}$ is a
Riemannian submanifold isometric to ${\bf R}^{n-1}$ for each $t\in
I$. By the equality $\al\circ \g(\cdot,\,t)=\g(\cdot,\,t)\circ
\al$ for every $\al\in G$, $B_t$ is exactly an orbit of the $G$
action, so it is a submanifold of $M$. Remeber that the map
$\g(\cdot,\,t):B\to B_t$ is surjective and $G$-equivariant.  By
the claim in the preceding paragraph, this map is one-to-one.
Hence, it gives a diffeomorphism of $B$ onto $B_t$. Since $G$ acts
effectively and isometrically on $B_t$, the claim follows from
Fact 1.3.

The left part of the proof is similarly to that of Lemma 3.2.
There also holds that that for each $(p,\,t)\in B\times I$
\[d\,(B,\,B_t)=d\,\bigl(\g(p,\,t),\,B\bigr)=d\,\bigl(\g(p,\,t),\,p\bigr)=t.\]
And the geodesic $\{\g(p,\,t):t\in I\}$ is perpendicular to $B_t$ at point $\g(p,\,t)$.\hf \\

The proof of Theorem 1.3 follows from Lemma 5.2. \hf \\

%%%%%%%%%%%%%%%%%%%%%%%% Acknowledgement %%%%%%%%%%%%%%%%%

\nd {\bf Acknowledgements}

\nd The second named author and the third are supported in part by
the National Natural Science Foundation of China No. 10671096 and
No. 10601053, respectively. The third named author would like to
express his sincere gratitude to Beijing International
Mathematical Research Center for her hospitality and financial
support during the course of this work.

%%%%%%%%%%%%%%%%%%% References %%%%%%%%%%%%%%%%%%%%%%%%%%%%%%%

\hspace{-0.70cm} {\sc Department of Mathematics}

\nd{\sc University of Science and Technology of China}

\nd {\sc Hefei 230026 China}

\nd {\sc E-mail addresses}: {\sc Zhi Chen} ({\tt
zzzchen@ustc.edu.cn}), {\sc Yiqian Shi} ({\tt yqshi@ustc.edu.cn}),
{\sc Bin Xu} ({\tt bxu@ustc.edu.cn}).

\end{document}